\documentclass[11pt]{amsart}
\usepackage{graphicx, amscd, amsfonts, amsthm}
\title{Spaces of embeddings of compact polyhedra into 2-manifolds}
\author{Tatsuhiko Yagasaki}
\address{Department of Mathematics, Kyoto Institute of Technology, Matsugasaki, Sakyoku, Kyoto 606, Japan}
\email{yagasaki@ipc.kit.ac.jp}

\newtheorem{theorem}{Theorem}[section]
\newtheorem{proposition}{Proposition}[section] 
\newtheorem{corollary}{Corollary}[section] 
\newtheorem{lemma}{Lemma}[section]
\newtheorem*{claim}{Claim}

\theoremstyle{definition}
\newtheorem{remark}{Remark}[section]

\newtheorem{fact}{Fact}[section]

\def \cal {\mathcal}
\def \phi {\varphi}

\begin{document}
\baselineskip 6 mm

\thispagestyle{empty}

\maketitle
\begin{abstract}
Let $M$ be a PL 2-manifold and $X$ be a compact subpolyhedron of $M$ and let ${\cal E}(X, M)$ denote the space of embeddings of $X$ into $M$ with the compact-open topology. In this paper we study an extension property of embeddings of $X$ into $M$ and show that the restriction map from the homeomorphism group of $M$ to ${\cal E}(X, M)$ is a principal bundle. As an application we show that if $M$ is a Euclidean PL $2$-manifold and $\dim X \geq 1$ then the triple $({\mathcal E}(X,M)$, ${\mathcal E}^{{\rm LIP}}(X,M)$, ${\mathcal E}^{{\rm PL}}(X, M))$ is an $(s,\Sigma,\sigma)$-manifold, where ${\mathcal E}_K^{{\rm LIP}}(X,M)$ and ${\mathcal E}_K^{{\rm PL}}(X, M)$ denote the subspaces of Lipschitz and PL embeddings.
\end{abstract} 

\section{Introduction}

The investigation of the topology of the homeomorphism groups of compact 2-manifolds \cite{Ha, HD, LM} included the use of conformal mappings in order to develop some extension properties of embeddings of a circle into an annulus and proper embeddings of an arc into a disk. In this paper we establish a similar extension property of embeddings of trees into a disk. Since every graph can be decomposed into ads (cones over finite points) and arcs connecting them, this implies an extension property of embeddings of compact polyhedra into 2-manifolds. 

Suppose $M$ is a PL 2-manifold and $K \subset X$ are compact subpolyhedra of $M$. Let ${\cal E}_K(X, M)$ denote the space of embeddings $f : X \hookrightarrow M$ with $f|_K = id$, equipped with the compact-open topology. An embedding $f : X \hookrightarrow M$ is said to be proper if $f(X \cap \partial M) \subset \partial M$ and $f(X \cap {\rm Int} \,M) \subset {\rm Int} \, M$. Let ${\cal E}_K(X, M)^{\ast}$ denote the subspace of proper embeddings in ${\cal E}_K(X, M)$, and let ${\cal E}_K(X, M)^{\ast}_0$ denote the connected component of the inclusion $i_X : X \subset M$ in ${\cal E}_K(X, M)^{\ast}$. Our result is summarized in the next statement. 

\begin{theorem}
For every $f \in {\mathcal E}_K(X, M)^{\ast}$ and every neighborhood $U$ of $f(X)$ in $M$, there exist a neighborhood ${\mathcal U}$ of $f$ in ${\mathcal E}_K(X, M)^{\ast}$ and a map $\phi : {\mathcal U} \to {\mathcal H}_{K \cup (M \setminus U)}(M)_0$ such that $\phi(g)f = g$ for each $g \in {\mathcal U}$ and $\phi(f) = id_M$.
\end{theorem}

Let ${\cal H}_X(M)$ denote the group of homeomorphisms $h$ of $M$ onto itself with $h|_X =id$, equipped with the compact-open topology. Let ${\cal H}(M)_0$ denote the identity component of ${\cal H}(M)$. In the study of the homotopy type of ${\cal H}_X(M)_0$ and ${\cal E}_K(X, M)_0$ the restriction map $\pi : {\cal H}_K(M)_0 \to {\cal E}_K(X, M)^\ast_0$ plays an important role (cf. \cite{Ep}). The above extension maps yield local sections of this restriction map. 

\begin{corollary}
For any open neighborhood $U$ of $X$ in $M$, the restriction map $\pi : {\cal H}_{K \cup (M \setminus U)}(M)_0 \to {\cal E}_K(X, U)_0^{\ast}$, $\pi(f) = f|_X$, is a principal bundle with the fiber ${\cal G} \equiv {\cal H}_{K \cup (M \setminus U)}(M)_0 \cap {\mathcal H}_X(M)$, where the subgroup ${\mathcal G}$ acts on ${\cal H}_{K \cup (M \setminus U)}(M)_0$ by right composition. 
\end{corollary}

As an application of Extension Theorem 1.1 we can study the embedding space ${\cal E}_K(X, M)$ from the viewpoint of infinite dimensional topology (see \S4 for basic terminologies). In \cite{SW} K.\,Sakai and R.\,Y.\,Wong showed the $(s, \Sigma, \sigma)$-stability property of triples of spaces of embeddings of compact polyhedra and subspaces of Lipschitz and PL embeddings, and posed the question whether these triples are $(s, \Sigma, \sigma)$-manifolds. The 1-dimensional case is studied in \cite{SK}. In this paper we will consider the 2-dimensional case and answer the question affirmatively. 

Let ${\cal E}_K^{{\rm PL}}(X,M)$ denote the subspace of PL-embeddings. When $M$ is a Euclidean PL 2-manifold, let ${\cal E}_K^{{\rm LIP}}(X,M)$ denote the subspace of Lipschitz embeddings. The Extension Theorem enables us to reduce the ANR-property and the homotopy negligibility of embedding spaces to the ones of the homeomorphism groups. Using the characterization of $(s, \Sigma, \sigma)$-manifold \cite{Ya} we have the following result. 

\begin{theorem}
Suppose $M$ is a Euclidean PL $2$-manifold and $K \subset X$ are compact subpolyhedra of $M$. If $\dim (X \setminus K) \geq 1$, then the triple $({\mathcal E}_K(X,M)$, ${\mathcal E}_K^{{\rm LIP}}(X,M)$, ${\mathcal E}_K^{{\rm PL}}(X, M))$ is an $(s,\Sigma,\sigma)$-manifold. 
\end{theorem}

Further applications of Corollary 1.1 to the study of ${\cal H}_X(M)$ and ${\cal E}_K(X, M)$ will be given in a succeeding paper. We conclude this section with some remarks. In Section 2 we study the extension property of embeddings of a tree into a disk. Section 3 contains the proofs of Theorem 1.1 and Corollary 1.1. The final section 4 contains the proof of Theorem 1.2. Throughout the paper spaces are assumed to be separable and metrizable. A Euclidean PL $n$-manifold is a subpolyhedron of some Euclidean space ${\Bbb R}^m$ which is a PL-manifold with respect to the induced triangulation and is equipped with the metric induced from the standard metric of ${\Bbb R}^m$. When $M$ is an orientable manifold, ${\cal H}_+(M)$  denote the subspace of orientation preserving homeomorphisms of $M$. Finally $i_X : X \subset Y$ denotes the inclusion map.

\section{Extension property of embeddings of trees into disks} 
In this section we will study some extension properties of embeddings of trees into disks. The proper embedding case is a consequence of a direct application of the conformal mapping theorem on simply connected domains (cf.\,\cite{LM}). Thus our interest is in the case of embeddings into the interior of a disk, where we need to apply the conformal mapping theorem on a doubly connected domain one boundary circle of which is collapsed to a tree. 

Throughout the section we will work on the plane ${\Bbb C}$ ($= {\Bbb R}^2)$ and use the following notations: For $z \in {\Bbb C}$ and $r > 0$, $D(z, r) = \{ x \in {\Bbb C} \, : \, |z - x| \leq r \}$, $O(z, r) = \{ x \in {\Bbb C} \, : \, |z - x| < r \}$, $C(z, r) = \{ x \in {\Bbb C} \, : \, |z - x| = r \}$, and $D(r) = D(0, r)$, $O(r) = O(0, r)$, $C(r) = C(0, r)$. For $0 < r < s$, $A(r, s) = \{ x \in {\Bbb C} \, : \, r \leq |x| \leq s \}$. For $A \subset {\Bbb C}$ and $\varepsilon > 0$, $O(A, \varepsilon) = \{ x \in {\Bbb C} \, : \, |x - y| < \varepsilon$ for some $y \in A \}$ (the $\varepsilon$-neighborhood of $A$).  

\subsection{Proper embeddings of trees into a disk} \mbox{} 

First we recall the conformal mapping theorem on simply connected domains normalized by the three points boundary condition. Consider the family ${\cal J} = \{ (J, w_1, w_2, w_3) \, : \, J$ is a simple closed curve in ${\Bbb C}$ and $w_1, w_2, w_3 \in J$ are three distinct points lying on $J$ in counterclockwise order (with respect to the orientation induced from ${\Bbb C}).\}$ A sequence $\{ A_n \}_{n \geq 1}$ of subsets of ${\Bbb C}$ is said to be uniformly locally connected if for each $\varepsilon > 0$ there exists a $\delta > 0$ such that for any $n \geq 1$ and any $x$, $y \in A_n$ with $|x - y| < \delta$ there exists an arc $\alpha$ in $A_n$ with connecting $x$ and $y$ and ${\rm diam} \, \alpha < \varepsilon$.

\begin{fact}
Let $z_1, z_2, z_3 \in C(1)$ be the fixed three points lying on $C(1)$ in counterclockwise order.\\
(i) (\cite[Corollary 2.7]{Po}) For every $(J, w_1, w_2, w_3) \in {\cal J}$ there exists a unique $\phi = \phi(J, w_1, w_2, w_3) \in {\cal E}(D(1), {\Bbb C})$ such that $\phi$ maps $O(1)$ conformally onto the interior of $J$, $\phi(C(1)) = J$ and $\phi(z_i) = w_i \ (i = 1, 2, 3)$. \\ (ii) If a sequence $(J_n, w_1(n), w_2(n), w_3(n)) \ (n \geq 1)$ converges to $(J, w_1, w_2, w_3)$ in the following sense, then $\phi(J_n, w_1(n), w_2(n), w_3(n))$ converges uniformly to $\phi(J, w_1, w_2, w_3)$: 

$(\ast)$ $J_n$ converges to $J$ with respect to the Hausdorff metric, $\{ J_n \}$ is uniformly locally connected, and $w_i(n) \to w_i \ (i = 1, 2, 3)$.
\end{fact}

For the statement (ii) we refer to the proof of \cite[Theorem 2.1, Proposition 2.3]{Po} (also see the proof of Lemma 2.3). 

\begin{lemma} Suppose $D$ is a disk and $C = \partial D$. \\ {\rm (i)} {\rm (cf. \cite[Lemma 3]{LM})} There exists a map $\Phi : {\cal E}(C, {\Bbb C}) \to {\cal E}(D, {\Bbb C})$ such that $\Phi(f)|_C = f \ (f \in {\cal E}(C, {\Bbb C}))$. \\ 
{\rm (ii)} {\rm (cf. \cite[Lemma 5]{LM})} Suppose $T$ is a tree embedded into a disk $D$ such that $T \cap C$ coincides with the set of terminal vertices of $T$. Then there exists a map $\Psi : {\cal E}_{T \cap C}(T, D)^\ast \to {\cal H}_{\partial}(D)$ such that $\Psi(f)|_{T} = f \ (f \in {\cal E}_{T \cap C}(T, D)^\ast)$ and $\Psi(i_T) = id_D$.
\end{lemma}

\begin{proof}
We may assume that $D = D(1)$. Let $z_1, z_2, z_3 \in C(1)$ be as in Fact 2.1. 

(i) Let ${\cal E}^{\pm} = \{ f \in {\cal E}(C(1), {\Bbb C}) \, : \, f$ preserves (reverses) orientation$\}$. If $f \in {\cal E}^+(C(1), {\Bbb C})$, then $(f(C(1)), f(z_1), f(z_2), f(z_3)) \in {\cal J}$ and by Fact 2.1 we obtain $\phi(f) = \phi(f(C(1)), f(z_1), f(z_2), f(z_3)) \in {\cal E}(D(1), {\Bbb C})$. If $f_n \to f$ in ${\cal E}^+$, then $(f(C(1)), f(z_1), f(z_2), f(z_3))$ converges to $(f(C(1)), f(z_1), f(z_2), f(z_3))$ in the sense $(\ast)$ of Fact 2.1.(ii). Hence the map $\phi : {\cal E}^+ \to {\cal E}(D(1), {\Bbb C})$ is continuous. Let $c : {\cal H}(C(1)) \to {\cal H}(D(1))$ be the cone extension map and
let $\gamma : {\Bbb C} \to {\Bbb C}$ be the reflection $\gamma(z) = \overline{z}$. Then the extension map $\Phi$ is defined by $\Phi(f) = \phi(f)c(\phi(f)^{-1} f)$ for $f \in {\cal E}^+$ and $\Phi(f) = \gamma \Phi(\gamma f)$ for $f \in {\cal E}^-$.

(ii) The tree $T$ separates the disk $D(1)$ into subdisks $D_i$. By (i) each disk $D_i$ admits an extension map $\psi_i : {\cal E}(\partial D_i, {\Bbb C}) \to {\cal E}(D_i, {\Bbb C})$. Every $f \in {\cal E}_{T \cap C(1)}(T, D(1))^\ast$ can be extended to $\overline{f} \in {\cal E}_{C(1)}(T \cup C(1), D(1))$. The required extension map $\Psi$ is defined by $\Psi(f)|_{D_i} = \psi_i(\overline{f}|_{\partial D_i})$. To achieve $\Psi(i_T) = id_D$, replace $\Psi(f)$ by $\Psi(f)\Psi(i_T)^{-1}$. \end{proof}

In the proof of Theorem 1.1 we will apply the statement (ii) to the case where $T$ is an arc.

\subsection{Embeddings of trees into the interior of a disk} \mbox{} 

Suppose $T$ is a finite tree ($\neq$ 1 pt) embedded into $O(2)$. We will use the following notation: For $a$, $b \in T$, let $E_T(a, b)$ denote the unique arc in $T$ connecting $a$ and $b$. Let $\{ v_1, \cdots, v_n \}$ be the collection of end vertices of $T$. We can choose disjoint arcs $\alpha_1, \cdots, \alpha_n$ in $D(2)$ such that each $\alpha_i$ connects $v_i$ with a point $a_i$ in $C(2)$ and ${\rm Int} \, \alpha_i \subset O(2) \setminus T$. We can arrange the ordering of $v_i$'s so that $a_1, \cdots, a_n$ lie on $C(2)$ in counterclockwise order. The labeling is unique up to the cyclic permutations. Note that $T$ does not meet the interior of the disk surrounded by the simple closed curve $\alpha_i \cup E_T(v_i, v_{i+1}) \cup \alpha_{i+1} \cup a_ia_{i+1}$, where $v_{n+1} = v_1$ and $a_{n+1} = a_1$. 

\begin{lemma} $($\cite[Ch.V, \S1, Theorems 1.1, 1.2]{Go}$)$ There exists a unique real number $r$, $0 < r < 2$, and a unique map $h : A(r, 2) \to D(2)$ such that $h : {\rm Int} \, A(r, 2) \to O(2) \setminus T$ is a conformal map and $h(2) = 2$. Furthermore, the map $h$ satisfies the following conditions : {\rm (i)} $h$ maps $C(2)$ homeomorphically onto $C(2)$,
{\rm (ii)} $h(C(r)) = T$ and there exists a unique collection of points $\{ u_1, \cdots, u_n \}$ lying on $C(r)$ in counterclockwise order such that $h$ maps each circular arc $u_iu_{i+1}$ homeomorphically onto the arc $E_T(v_i, v_{i+1})$, where $u_{n+1} = u_1$.
\end{lemma}

We refer to \cite[Ch2. Theorem 2.1]{Po} for the extension to boundary and \cite[Ch2. \S1 Prime End Theorem, \S\S4, 5]{Po} and \cite[p.40]{Go} for the correspondence between prime ends and boundary points. 
Let ${\cal E} = {\cal E}(T, O(2))$. For each $f \in {\cal E}$ the image $f(T)$ is a tree in $O(2)$. Hence by Lemma 2.2 there exists a unque real number $r_f$, $0 < r_f < 2$, and a unique map $h_f : A(r_f, 2) \to D(2)$ such that $h_f : {\rm Int} \, A(r_f, 2) \to O(2) \setminus f(T)$ is a conformal map and $h_f(2) = 2$. For $0 < r < 2 $ let $\phi_r : A(1, 2) \to A(r, 2)$ denote the radial map defined by $\phi_r(x) = ((2 - r)(|x| - 1) + r)x/|x|$, and let ${\cal C}(A(1, 2), D(2))$ denote the space of continuous maps from $A(1, 2)$ to $D(2)$, with the compact-open topology. We have $h_f \phi_{r_f} \in {\cal C}(A(1, 2), D(2))$. 

\begin{lemma}
The map $\Psi : {\cal E}(T, O(2)) \to {\Bbb R} \times {\cal C}(A(1, 2), D(2))$, $\Psi(f) = (r_f, h_f \phi_{r_f})$, is continuous.
\end{lemma}

This continuity property can be verified using the length distortion under conformal mapping \cite[Proposition 2.2]{Po}. When $L$ is a rectifiable (possibly open) curve in ${\Bbb R}^2$, we denote the length of $L$ by $\Lambda(L)$. 

\begin{proof}
Suppose $f_n \to f$ in ${\cal E}$. It suffices to show that the sequence $(r_n, h_n\phi_{r_n}) \equiv (r_{f_n}, h_{f_n}\phi_{r_{f_n}})$ has a subsequence $(r_{n_k}, h_{n_k}\phi_{r_{n_k}})$ such that $r_{n_k} \to r_f$ and $h_{n_k} \phi_{r_{n_k}}$ converges uniformly to $h_f \phi_{r_f} $. 

Let $R_0 > 2 \ (=$ the radius of $D(2)$) and $\varepsilon(\rho) = 2 \pi R_0/ \sqrt{\log \, (1/\rho)} \ (0 < \rho < 1)$. 

(i) Passing to a subsequence we may assume $r_n \to r_0$ for some $r_0, 0 \leq r_0 \leq 2$.
First we will show that $0 < r_0 < 2$. (a) Suppose $r_0 = 2$. 
Take $\rho$, $0 < \rho < 1$, such that $\varepsilon(\rho) < d(f(T), C(2))$. Choose $n \geq 1$ such that $\varepsilon(\rho) < d(f_n(T), C(2))$ and $|r_n - r_0| < \rho$.
We can apply \cite[Proposition 2.2]{Po} for any point $c \in C(2)$ (with $R = 2$) to find $\rho_0, \rho < \rho_0 < \sqrt{\rho}$, such that $\Lambda(h_n(L)) < \varepsilon(\rho)$, where $L$ is one of the two arc components of $C(c, \rho_0) \cap A(r_n, 2)$ which connects $C(r_n)$ and $C(2)$.
This implies $d(f_n(T), C(2)) < \varepsilon(\rho)$, a contradiction. 
(b) Suppose $r_0 = 0$. Take $\rho$, $0 < \rho < 1$, such that $\varepsilon(\rho) < {\rm diam} \, f(T)$. Choose $n \geq 1$ such that $\varepsilon(\rho) < {\rm diam} \, f_n(T)$ and $r_n < \rho$. By \cite[Proposition 2.2]{Po} there exists $\rho_0$, $\rho < \rho_0 < \sqrt{\rho}$ such that $\Lambda(h_n(C(\rho_0))) < \varepsilon(\rho)$. Since $f_n(T)$ is contained in the interior of the circle $h_n(C(\rho_0))$, we have ${\rm diam} \, f_n(T) < \varepsilon(\rho)$, a contradiction.

(ii) Next we will show that the sequence $h_n : A(r_n, 2) \to D(2) \ (n \geq 1)$ is equicontinuous, that is, for every $\varepsilon > 0$ there exists a $\rho > 0$ such that $|h_n(z) - h_n(w)| < \varepsilon$ for any $n \geq 1$ and $z$, $w \in A(r_n, 2)$ with $|z - w| < \rho$. Let $\varepsilon > 0$ be given.
We may assume that $\varepsilon < d(C(2), f_n(T))$ for each $n \geq 1$. Since the sequence $C(2)$, $f_n(T) \ (n \geq 1)$ is uniformly locally connected, there exists a $\delta$, $0 < \delta < \varepsilon/2$, such that if $z$, $w \in f_n(T)$ (respectively $C(2)$) and $|z - w| < \delta$, then there exists an arc $A$ in $f_n(T)$ (respectively $C(2)$) connecting $z$ and $w$ and with ${\rm diam} \, A < \varepsilon/2$.
Choose $\rho$, $0 < \rho < 1$, such that $\varepsilon(\rho) < \delta$ and $2\sqrt{\rho} < 2 - \max_{n \geq 0} \, r_n$. Suppose $z$, $w \in A(r_n, 2)$ and $|z - w| < \rho$. By \cite[Proposition 2.2]{Po} (with $c = z$) we have $\rho_0$, $\rho < \rho_0 < \sqrt{\rho}$, such that $\Lambda(h_n(L)) < \varepsilon(\rho)$, where $L = C(z, \rho_0) \cap A(r_n, 2)$. Since $z$, $w \in D \equiv D(z, \rho_0) \cap A(r_n, 2)$, it suffices to show that $\text{diam} \, h_n(D) < \varepsilon$. By the choice of $\rho$, $D(z, \rho_0)$ meet at most one of $C(2)$ and $C(r_n)$. If $D(z, \rho_0) \subset A(r_n, 2)$ or $D(z, \rho_0) \supset D(0, r_n)$, then $L = C(z, \rho_0)$ and $h_n(D)$ is a disk bounded by $h_n(L)$, so $\text{diam} \, h_n(D) < \varepsilon(\rho)$. Otherwise, $L$ is an arc connecting two points $P$, $Q$ with either (a) $P, Q \in C(2)$ or (b) $P, Q \in C(r_n)$. In both cases $|h_n(P) - h_n(Q)| \leq \Lambda(h_n(L)) < \delta$, hence by the choice of $\delta$, we have an arc $A$ in $C(2)$ (resp. $f_n(T)$) connecting $h_n(P)$ and $h_n(Q)$ and $\text{diam} \, A < \varepsilon/2$. In the case (a) $h_n(L)$ separates $D(2)$ into the subdisk $h_n(D)$ and another subdisk. Since $h_n(D) \cap f_n(T) = \emptyset$ and $d(C(2), f_n(T)) > \varepsilon$, the Jordan curve $h_n(L) \cup A$ bounds the disk $h_n(D)$, so $\text{diam} \, h_n(D) < \varepsilon$. In the case (b) the Jordan curve $h_n(L) \cup A$ bounds a disk $E$ in $D(2)$ with $\text{diam} \, E < \varepsilon$. Since $h_n(A(r_n, 2) \setminus (D \cup C(r_n)))$ is contained in the exterior of $E$ and $h_n(\text{Int} \, D) \cap \partial E = \emptyset$, it follows that $h_n(\text{Int} \, D) = \text{Int} \, E \setminus f_n(T)$ so $h_n(D) = E$.

(iii) Since the sup-metric $d(\phi_{r_n}, \phi_{r_0}) = |r_n - r_0| \to 0 \ (n \to \infty)$, the sequence $h_n\phi_{r_n} \ (n \geq 1)$ is also equicontinuous.
By the Ascoli-Arzel\`a theorem, passing to a subsequence, we may assume that $h_n\phi_{r_n}$ converges to a map $h_0' : A(1, 2) \to D(2)$. Set $h_0 = h_0'\phi_{r_0}^{-1}$.
Then $h_0(A(r_0, 2)) = D(2)$, $h_0(C(2)) = C(2)$ $h_0(C(r_0)) = f(T)$ and $h_0(2) = 2$. Since the sequence of univalent analytic maps $h_n : \text{Int} \, A(r_n, 2) \to {\Bbb C}$ converges weakly uniformly to the map $h_0 : \text{Int} \, A(r_0, 2) \to {\Bbb C}$ (i.e., for each compact subset $K$ of ${\rm Int} \, A(r_0, 2)$, $h_n|_K$ ($n$ large) converges uniformly to $h_0|_K$) and $h_0$ is not constant, $h_0 : \text{Int} \, A(r_0, 2) \to {\Bbb C}$ is also a univalent analytic map \cite[Ch.3, Theorem 3.3]{Ve}.
It follows that $h_0(\text{Int} \, A(r_0, 2)) = O(2) \setminus f(T)$ and $h_0 : \text{Int} \, A(r_0, 2) \to O(2) \setminus f(T)$ is a conformal map, so $(r_0, h_0) = (r_f, h_f)$ by the uniqueness in Lemma 2.2. This completes the proof.
\end{proof}

Let $i : T \hookrightarrow O(2)$ denote the inclusion and set ${\cal E}_+ \equiv {\cal E}_+(T, O(2)) = \{ f \in {\cal E} \, :$ there exists an $h \in {\cal H}_+(D(2))$ with $hi = f \}$, which is an open neighborhood of $i$ in ${\cal E}$. 

\begin{proposition}
{\rm (i)} There exists a canonical map $\Phi = \Phi_T : {\mathcal E}_+ \to {\cal H}_+(D(2))$ such that $\Phi(f)i = f \ (f \in {\mathcal E}_+)$ and $\Phi(i) = id$.

{\rm (ii)} There exists a neighborhood ${\mathcal U}$ of \hspace{1pt} $i$ in ${\mathcal E}$ and a map $\phi : {\mathcal U} \to {\cal H}_{\partial}(D(2))$ such that $\phi(f)i = f \ (f \in {\mathcal U})$ and $\phi(i) = id_D$.
\end{proposition}

\begin{proof}
(i) Let $f \in {\mathcal E}_+$.
Comparing two maps $h_f \phi_{r_f}$, $fh_i \phi_{r_i} : C(1) \to f(T)$, we obtain a unique map $\Theta_0(f) \in {\mathcal H}_+(C(1))$ such that $h_f \phi_{r_f} \Theta_0(f) = fh_i \phi_{r_i}$. Extend $\Theta_0(f)$ radially to $\Theta(f) \in {\mathcal H}_+(A(1, 2))$ by $\Theta(f)(rz) = r\Theta_0(f)(z) \ (z \in C(1), 1 \leq r \leq 2)$. The required map $\Phi(f)$ is defined as the unique map $\Phi(f) \in {\cal H}_+(D(2))$ with $h_f \phi_{r_f}\Theta(f) = \Phi(f) h_i \phi_{r_i}$. In Claim below we will show that the map $\Theta_0$ is continuous. This implies the continuity of the map $\Phi$. 

(ii) Since $\Phi(i) = id$, if we take a sufficiently small neighborhood ${\mathcal U}$ of $i$, then $\Phi(f)|_{C(2)}$ is close to $id_{C(2)}$ for $f \in {\mathcal U}$, and we can use a collar of $C(2)$ in $D(2)$ and a local contraction of a neighborhood of $id_{C(2)}$ in ${\mathcal H}(C(2))$ to modify the map $\Phi|_{{\mathcal U}}$ to obtain the desired map $\phi$. \end{proof}

\begin{claim} 
The map $\Theta_0 : {\cal E}_+ \to {\cal H}_+(C(1))$ is continuous. 
\end{claim} 

\begin{proof}
Under the notations of Lemma 2.2, let $g_f = h_f \phi_{r_f}$ and $x_j(f) = \phi_{r_f}^{-1}(u_j)$. For the inclusion $i : T \subset D(2)$, we abbreviate as $g = g_i$ and $x_j = x_j(i)$. Let $L_j = x_j x_{j+1}$ (the circular arc in $C(1)$). Also let $\tilde{f} = \Theta_0(f)$. Note that $g_f$ is continuous in $f$ (Lemma 2.3), $g_f \tilde{f} = gf$, $\tilde{f}(x_j) = x_j(f) = g_f^{-1}(f(v_j))$ and that $g_f$ maps $\tilde{f}(L_j)$ homeomorphically onto $f(E_T(v_j, v_{j+1}))$. 

(1) First we will show the following statement: 

($\ast$) Suppose $f \in {\cal E}_+, U$ is any open neighborhood of $x_j(f)$ in ${\Bbb C}$ and $A_j$ is a small compact neighborhood of $x_j$ in $C(1)$ such that $g_f(\tilde{f}(A_j)) \cap g_f(A(1, 2) \setminus U) = \emptyset$ (hence $\tilde{f}(A_j) \subset U$). If $f'$ is sufficiently close to $f$, then $\tilde{f'}(A_j) \subset U$. In particular, $x_j(f) \in C(1)$ is continuous in $f$.

In fact, there exists an $\varepsilon > 0$ such that $O(fg(A_j), \varepsilon) \cap O(g_f(A(1, 2) \setminus U), \varepsilon) = \emptyset$. If $f'$ is sufficiently close to $f$ then the sup-metric $d(f', f) < \varepsilon$ and $d(g_{f'}, g_f) < \varepsilon$. Hence, $f'g(A_j) = g_{f'} \tilde{f'}(A_j)$ does not meet $g_{f'}(A(1, 2) \setminus U)$, so $g_{f'}\tilde{f'}(A_j) \subset U$. 

(2) To show that $\tilde{f}$ is continuous in $f$, let $f \in {\cal E}_+$ and $\varepsilon > 0$ be given. It suffices to show that for each $j = 1, \cdots, n$ there exists a small neighborhood ${\cal U}$ of $f$ in ${\cal E}_+$ such that $\tilde{f}$ and $\tilde{f'}$ are $\varepsilon$-close on $L_j$ for every $f' \in {\cal U}$. 

Set $U_j = O(x_j(f), \varepsilon/2)$ and $U_{j+1} = O(x_{j+1}(f), \varepsilon/2)$, and let $A_j$ and $A_{j+1}$ be small circular arc neighborhoods of $x_j$ and $x_{j+1}$ in $C(1)$ as in (1) with respect to $U_j$ and $U_{j+1}$ respectively. Set $K_j = cl(L_j \setminus (A_j \cup A_{j+1}))$ and choose small circular arc neighborhoods $C_j$ and $C_{j+1}$ of $x_j(f)$ and $x_{j+1}(f)$ in $C(1)$ such that $g_f\tilde{f}(K_j)$ meets neither $g_f(C_j)$ nor $g_f(C_{j+1})$. Choose $\delta_1 > 0$ such that $O(g_f\tilde{f}(K_j), \delta_1)$ meets neither $O(g_f(C_j), \delta_1)$ nor $O(g_f(C_{j+1}), \delta_1)$. By the compactness argument there exists $\delta$, $0 < \delta < \delta_1$, such that for any $x \in L_j$, $g_f(\tilde{f}(L_j) \setminus O(\tilde{f}(x), \varepsilon)) \cap O(g_f\tilde{f}(x), 2 \delta) = \emptyset$.

By (1) there exists a neighborhood ${\cal U}$ of $f$ in ${\cal E}_+$ such that if $f' \in {\cal U}$, then $\tilde{f'}(A_j) \subset U_j$, $\tilde{f'}(A_{j+1}) \subset U_{j+1}$, $\tilde{f'}(x_j) \in C_j$, $\tilde{f'}(x_{j+1}) \in C_{j+1}$ and $d(f, f') < \delta$, $d(g_{f'}, g_f) < \delta$. Since $\tilde{f'}$ is orientation preserving, $\tilde{f'}(x_j) \in C_j$ and $\tilde{f'}(x_{j+1}) \in C_{j+1}$, it follows that $\tilde{f'}(L_j) \subset \tilde{f}(L_j) \cup C_j \cup C_{j+1}$. If $x \in A_j$, then $\tilde{f'}(x), \tilde{f}(x) \in U_j$ so that $d(\tilde{f'}(x), \tilde{f}(x)) < \varepsilon$. For each $x \in A_{j+1}$ we have the same conclusion. Suppose $x \in K_j$. Since $g_{f'} \tilde{f'}(x) = f'g(x)$ is $\delta$-close to $fg(x) = g_f \tilde{f}(x) \in g_f \tilde{f}(K_j)$ and $g_{f'}(C_j) \subset O(g_f(C_j), \delta)$, we have $\tilde{f'}(x) \not\in C_j$. Similarly $\tilde{f'}(x) \not\in C_{j+1}$, and so $\tilde{f'}(x) \in \tilde{f}(L_j)$. Since $g_f\tilde{f}(x) = fg(x)$ is $\delta$-close to $f'g(x) = g_{f'} \tilde{f'}(x)$ and the latter is also $\delta$-close to $g_f \tilde{f'}(x)$, we have $g_f \tilde{f'}(x) \in O(g_f \tilde{f}(x), 2 \delta)$. Hence by the choice of $\delta$, $\tilde{f'}(x) \in O(\tilde{f}(x), \varepsilon)$. This completes the proof. 
\end{proof}

Finally we will see a symmetry property of the map $\Phi_T$ in Proposition 2.1 (i).
For $z \in C(1)$ let $\theta_z : {\Bbb C} \to {\Bbb C}$ denote the rotation $\theta_z(w) = z \cdot w$ and let $\gamma : {\Bbb R}^2 \to {\Bbb R}^2$ be the reflection, $\gamma (x, y) = (x, -y)$. 

\begin{lemma}
{\rm (i)} $\Phi_T(\theta_z f) = \theta_z \Phi_T(f) \ (f \in {\cal E}_+, z \in C(1))$. \\
{\rm (ii)} $\Phi_{\gamma(T)}(\gamma f \gamma) = \gamma \Phi_T(f) \gamma \ (f \in {\cal E}_+)$. In particular, if $T$ is a segment in the $x$-axis, then $\Phi_T(\gamma f) = \gamma \Phi_T(f) \gamma \ (f \in {\cal E})$. \end{lemma}

\begin{proof}
(i) Let $f \in {\cal E}_+$, $z \in C(1)$ and let $w_0 \in C(2)$ be the unique point such that $\theta_z h_f \theta_z^{-1}(w_0) = 2$. Under Lemma 2.2, $(r_f, \theta_z h_f \theta_z^{-1} \theta_w)$ corresponds to $\theta_z f$, where $w = w_0/2$. Thus $\Theta(\theta_z f) = \theta_w^{-1} \theta_z \Theta(f)$ and the conclusion follows from 
\[ \Phi(\theta_z f)h_i \phi_{r_i} = (\theta_z h_f \theta_z^{-1} \theta_w) \phi_{r_f} \Theta(\theta_z f) = (\theta_z h_f \theta_z^{-1} \theta_w) \phi_{r_f} \theta_w^{-1} \theta_z \Theta(f) = \theta_z h_f \phi_{r_f} \Theta(f) = \theta_z \Phi(f) h_i \phi_{r_i}. \] (ii) Since $(r_i, \gamma h_i \gamma)$ corresponds to $\gamma(T)$ and $(r_f, \gamma h_f \gamma)$ corresponds to $\gamma f(T)$, it follows that $\Theta_{\gamma(T)}(\gamma f \gamma) = \gamma \Theta_T(f) \gamma$. The conclusion follows from 
\[ (\gamma \Phi(f) \gamma) (\gamma h_i \gamma \phi_{r_i}) = \gamma (\Phi(f) h_i \phi_{r_i}) \gamma = \gamma (h_f \phi_{r_f} \Theta(f)) \gamma = (\gamma h_f \gamma \phi_{r_f})(\gamma \Theta(f) \gamma) \]
\end{proof}

\section{Extension property of embeddings of compact polyhedra into 2-manifolds}

In this section we prove Theorem 1.1 and Corollary 1.1. First we consider the case where $M$ is compact. 

\begin{lemma}
Suppose $M$ is a compact PL $2$-manifold and $K \subset X$ are compact subpolyhedra of $M$. Then there exists an open neighborhood ${\cal U}$ of $i_X$ in ${\cal E}_K(X, M)^{\ast}$ and a map $\phi : {\cal U} \to {\cal H}_K(M)$ such that $\phi(f)|_X = f \ (f \in {\cal U})$ and $\phi(i_X) = id_M$.
\end{lemma} 

\begin{proof}
We may assume that $K = \emptyset$, since if $\phi$ satisfies the above condition in the case where $K = \emptyset$ then we have $\phi({\cal U} \cap {\cal E}_K(X, M)^{\ast}) \subset {\cal H}_K(M)$ for any $K \subset X$.

(1) The case when $\partial \, M = \emptyset$: We fix a triangulation of $X$ and let $S_k \ (k = 0, 1, 2)$ denote the set of $k$-simplices of this triangulation and $X^{(1)}$ denote the 1-skeleton of $X$.
For each $\sigma \in S_1$ with ends $v$, $w$ we choose two disjoint subarcs $\sigma_v$, $\sigma_w$ of $\sigma$ with $v \in \sigma_v$, $w \in \sigma_w$ and a subarc $e_{\sigma}$ of ${\rm Int} \, \sigma$ with ${\rm Int} \, e_{\sigma} \supset cl \, (\sigma \setminus (\sigma_v \cup \sigma_w))$. For each $v \in S_0$ set $T_v = \{ v \} \cup (\cup_{v \in \sigma \in S_1} \, \sigma_v)$, which is an ad or a single point. We choose two disjoint families of closed disks $\{ D_v \}_{v \in S_0}$ and $\{ E_\sigma \}_{\sigma \in S_1}$ in $M$ such that (i) $T_v \subset {\rm Int} \, D_v \ (v \in S_0)$ and (ii) $X^{(1)} \cap E_{\sigma} = e_{\sigma}$ and ${\rm Int} \, e_{\sigma} \subset {\rm Int} \, E_{\sigma}$ (i.e., $e_{\sigma}$ is a proper arc of $E_{\sigma})$. 

\[ \fbox{Figure 1.a} \]
\vskip 5mm 

By Proposition 2.1.(ii) for each $v \in S_0$ there exists a neighborhood ${\cal V}_v$ of $i_{T_v}$ in ${\cal E}(T_v, {\rm Int} \, D_v)$ and an extension map $\alpha_v : {\cal V}_v \to {\cal H}_{\partial}(D_v)$. In turn, by Lemma 2.1.(ii) for each $\sigma \in S_1$ there exists a neighborhood ${\cal W}_\sigma$ of $i_{e_\sigma}$ in ${\cal E}_{\partial e_{\sigma}}(e_{\sigma}, E_{\sigma})^{\ast}$ and an extension map $\beta_{\sigma} : {\cal W}_\sigma \to {\cal H}_{\partial}(E_{\sigma})$. If ${\cal U}$ is a sufficiently small neighborhood of $i_X$ in ${\cal E}(X, M)$, then for any $f \in {\cal U}$ we have $f|_{T_v} \in {\cal V}_v$ for every $v \in S_0$ and we can define a map $\lambda : {\cal U} \to {\cal H}(M)$ by \[ \lambda(f) = \cases \alpha_v(f|_{T_v}) & \mbox{on} \ D_v, \\ id & \mbox{on} \ M \setminus \cup_v \, D_v. \endcases \] Since $\lambda(i_X) = id_M$ and $\lambda(f)^{-1}f|_{T_v} = i_{T_v} \ (v \in S_0)$, if ${\cal U}$ is small enough, then $\lambda(f)^{-1}f$ is sufficiently close to $i_X$ so that $\lambda(f)^{-1}f|_{e_\sigma} \in {\cal W}_\sigma$.
Hence we can define a map $\mu : {\cal U} \to {\cal H}(M)$ by \[ \mu(f) = \cases
\beta_{\sigma}(\lambda(f)^{-1}f|_{e_\sigma}) & {\rm on} \ E_\sigma, \\ id & {\rm on} \ M \setminus \cup_\sigma \, E_\sigma. \endcases \]
Then $\mu(i_X) = id_M$ and $\hat{f} \equiv \mu(f)^{-1}\lambda(f)^{-1}f$ is equal to the identity map on $X^{(1)}$ for each $f \in {\cal U}$. Since $\hat{f}(\sigma) = \sigma \ (\sigma \in S_2)$, we can define a map $\nu : {\cal U} \to {\cal H}(M)$ by
$\nu(f)|_X = \hat{f}$ and $\nu(f)|_{M \setminus X} = id$. Since $\nu(i_X) = id_M$ and $\nu(f)^{-1}\mu(f)^{-1}\lambda(f)^{-1}f = i_X$, the map $\phi : {\cal U} \to {\cal H}(M)$, $\phi(f) = \lambda(f)\mu(f)\nu(f) \ (f \in {\cal U})$ satisfies the desired conditions. 

\[ \fbox{Figure 1-b.}  \hspace{2cm} \fbox{Figure 1-c.} \]
\vskip 5mm 

(2) The case when $\partial \, M \neq \emptyset$: We can use the double $N = M \cup_{\partial M} M$. Since $X$ is a subpolyhedron of $M$, $Y = X \cap \, \partial M$ is also a subpolyhedron of $\partial M$.

(i) By (1) (where $K \neq \emptyset$) we have a neighborhood ${\cal V}_0$ of $i_{X \cup \, \partial M}$ in ${\cal E}_{\partial M}(X \cup \partial M, N)$ and an extension map $\psi_0 : {\cal V}_0 \to {\cal H}_{\partial M}(N)$. We can extend every $f \in {\cal E}_Y(X, M)^\ast$ to an $f_0 \in {\cal E}_{\partial M}(X \cup \partial M, N)$ by the identity on $\partial M$. If ${\cal V}$ is a small neighborhood of $i_X$ in ${\cal E}_Y(X, M)^\ast$, then for every $f \in {\cal V}$ we have $f_0 \in {\cal V}_0$, so $\psi(f_0)$ is defined and $\psi_0(f_0)(M) = M$.
Thus we have an extension map $\psi : {\cal V} \to {\cal H}_{\partial M}(M)$, $\psi(f) = \psi_0(f_0)|_M$.

(ii) Since ${\cal H}(\partial M)$ is locally contractible, using a collar of $\partial M$ in $M$, we have a neighborhood ${\cal W}$ of $id_{\partial M}$ in ${\cal H}(\partial M)$ and a map $F : {\cal W} \to {\cal H}(M)$ such that $F(g)|_{\partial M} = g \ (g \in {\cal W})$ and $F(id_{\partial M}) = id_M$. We can easily verify a 1-dimensional version of Lemma 3.1 and find a neighborhood ${\cal W}_0$ of $i_Y$ in ${\cal E}(Y, \partial M)$ and an extension map $\lambda_0 : {\cal W}_0 \to {\cal H}(\partial M)$. We may assume that $\lambda_0({\cal W}_0) \subset {\cal W}$. 
Hence if ${\cal U}$ is a small neighborhood of $i_X$ in ${\cal E}(X, M)^\ast$, then we have a map $\lambda : {\cal U} \to {\cal H}(M)$, $\lambda(f) = F(\lambda_0(f|_Y))$.
Then $\lambda(f)|_Y = f|_Y \ (f \in {\cal U})$ and $\lambda(id_X) = id_M$. If ${\cal U}$ is small, then we have $\lambda(f)^{-1}f \in {\cal V}$ and the required extension map $\phi : {\cal U} \to {\cal H}(M)$ is defined by $\phi(f) = \lambda(f)\psi(\lambda(f)^{-1}f)$. 
\end{proof}

\begin{lemma}
If $M$ is a compact PL $2$-manifold and $X$ is a compact subpolyhedron of $M$, then ${\cal H}_X(M)$ is an ANR. 
\end{lemma}

\begin{proof}
Let $\pi : {\cal H}(M) \to {\cal E}(X, M)^\ast$, $\pi(h) = h|_X$, denote the restriction map. By Lemm 3.1 (with $K = \emptyset$) there exists an open neighborhood ${\cal U}$ of $i_X$ in ${\cal E}(X, M)^\ast$ and a map $\phi : {\cal U} \to {\cal H}(M)$ such that $\phi(f)|_X = f$. Then $\Phi : {\cal U} \times {\cal H}_X(M) \cong \pi^{-1}({\cal U})$, $\Phi(f, h) = \phi(f)h$, is a homeomorphism with the inverse $\Phi^{-1}(k) = (k|_X, \phi(k|_X)^{-1} k)$. Since ${\cal H}(M)$ is an ANR \cite{LM} and $\pi^{-1}({\cal U})$ is open in ${\cal H}(M)$, ${\cal H}_X(M)$ is also an ANR. 
\end{proof} 

\noindent {\bf Proof of Theorem 1.1.} Theorem 1.1 can be reduced to Lemma 3.1 by the following observations: \\
(i) Since there exists an $h \in {\cal H}_{K \cup (M \setminus U)}(M)_0$ such that $hf$ is a PL embedding (cf. \cite[Appendix]{Ep}) we may assume that $f$ is a PL-embedding. Replacing $X$ by $f(X)$, we may assume that $f = i_X : X \subset M$. \\ 
(ii) Taking a compact PL-submanifold neighborhood $N$ of $X$ in $U$ and replacing $(M, X, K)$ by $(N, X \cup \, {\rm Fr}_M N, K \cup \, {\rm Fr}_M N)$, we may assume that $M$ is compact and $U = M$. \\ 
(iii) If $M$ is compact then ${\cal H}_K(M)_0$ is open in ${\cal H}_K(M)$ by Lemma 3.2. Hence we can take a smaller ${\cal U}$ to attain $\phi({\cal U}) \subset {\cal H}_K(M)_0$. \qed 
\vskip 3mm 
\noindent {\bf Proof of Corollary 1.1.}
Let $f \in {\cal E}_K(X, U)^{\ast}_0$ and let ${\cal U}_f$, $\phi_f$ be as in Theorem 1.1. If ${\cal U}_f \cap {\rm Im} \, \pi \neq \emptyset$ then ${\cal U}_f \subset {\rm Im} \, \pi$. In fact, if $h \in {\cal H}_{K \cup (M \setminus U)}(M)_0$ and $\pi(h) = h|_X \in {\cal U}_f$, then for any $g \in {\cal U}_f$ we have $g = \pi(\phi_f(g) \phi_f(h|_X)^{-1}h)$. Hence ${\rm Im} \, \pi$ is clopen in ${\cal E}_K(X, U)^\ast_0$, so ${\rm Im} \, \pi = {\cal E}_K(X, U)^{\ast}_0$ and ${\cal U}_f \subset {\cal E}_K(X, U)^{\ast}_0$. Choose an $h_f \in {\cal H}_{K \cup (M \setminus U)}(M)_0$ with $h_f|_X = f$ and define a local trivialization $\Phi : {\cal U}_f \times {\cal G} \cong \pi^{-1}({\cal U}_f)$ by $\Phi(g, h) = \phi_f(g) h_f h$.  \qed 
\vskip 3mm
By a similar argument we can also show the following statements. 

\begin{corollary} Suppose $K \subset Y \subset X$ are compact subpolyhedra of a PL $2$-manifold $M$. \\ 
{\rm (i)} For any open neighborhood $U$ of $X$ in $M$ the restriction map $\pi : {\cal H}_{K \cup (M \setminus U)}(M) \to {\rm Im} \, \pi \subset {\cal E}_K(X, U)^\ast$ is a principal bundle with the fiber ${\cal H}_{X \cup (M \setminus U)}(M)$ and ${\rm Im} \, \pi$ is clopen in ${\cal E}_K(X, U)^\ast$. \\ 
{\rm (ii)} The restriction map $p : {\cal E}_K(X, M)^\ast \to {\rm Im} \, p \subset {\cal E}_K(Y, M)^\ast$ is locally trivial and ${\rm Im} \, p$ is clopen in ${\cal E}_K(Y, M)^\ast$.
\end{corollary} 

\section{The spaces of embeddings into 2-manifolds} 

In this final section we will prove Theorem 1.2. 

\subsection{Basic facts on infinite-dimensional manifolds} \mbox{} \vskip 1mm
First we recall some basic facts on infinite-dimensional manifolds. As for the model spaces we follow the standard convension: $s = (-\infty, \infty)^{\infty} \ (\cong \ell_2)$, $\Sigma = \{(x_n) \in s \, : \, \sup_n\,|x_n| < \infty \}$, $\sigma = \{(x_n) \in s \, : \, x_n = 0$ (almost all $n)\}$. 
A triple $(X, X_1, X_2)$ means a triple of a space $X$ and subspaces $X_1 \supset X_2$.
A triple $(X, X_1, X_2)$ is said to be a $(s, \Sigma, \sigma)$-manifold if each point of $X$ admits an open neighborhood $U$ in $X$ and an open set $V$ in $s$ such that $(U, U \cap X_1, U \cap X_2) \cong (V, V \cap \Sigma, V \cap \sigma)$ (a homeomorphism of triples). 
 In \cite{Ya} we have obtained a characterization of $(s, \Sigma, \sigma)$-manifolds in terms of some class conditions, a stability condition and the homotopy negligible complement condition. A space is $\sigma$-(fd-)compact if it is a countable union of (finite dimensional) compact subsets. A triple $(X, X_1, X_2)$ is said to be $(s, \Sigma, \sigma)$-stable if $(X \times s, X_1 \times \Sigma, X_2 \times \sigma) \cong (X, X_1, X_2)$. We say that a subset $Y$ of $X$ has the homotopy negligible (h.n.) complement in $X$ if there exists a homotopy $\phi_t : X \to X \ (0 \leq t \leq 1)$ such that $\phi_0 = id_X$ and $\phi_t(X) \subset Y \ (0 < t \leq 1)$. The homotopy $\phi_t$ is called an absorbing homotopy of $X$ into $Y$.  

\begin{fact}
(i) $Y$ has the h.n.\ complement in $X$ iff each point $x \in X$ has an open neighborhood $U$ and a homotopy $\phi : U \times [0, 1] \to X$ such that $\phi_0 = i_U : U \subset X$ and $\phi_t(U) \subset Y \ (0 < t \leq 1)$. \\
(ii) If $Y$ has the h.n.\ complement in $X$, then $X$ is an ANR iff $Y$ is an ANR by \cite{Han}. \\
(iii) (\cite{To}) Suppose $X$ is an ANR. Then $Y$ has the h.n.\ complement in $X$ iff for any open set $U$ of $X$ the inclusion $U \cap Y \subset U$ is a weak homotopy equivalence. Hence if both $Y \subset X$ and $Z \subset Y$ have the h.n.\ complement, then so does $Z \subset X$.
\end{fact}

In (i) $U \cap Y$ has the h.n.\ complement in $U$ and local absorbing homotopies can be uniformized to a global one \cite{Mi}. 

We will apply the following characterization of $(s, \Sigma, \sigma)$-manifolds \cite{Ya}. 

\begin{proposition}
A triple $(X, X_1, X_2)$ is an $(s, \Sigma, \sigma)$-manifold iff \\ 
{\rm (i)} $X$ is a separable completely metrizable ANR, $X_1$ is $\sigma$-compact and $X_2$ is $\sigma$-fd-compact, \\ 
{\rm (ii)} $X_2$ has the h.n.\ complement in $X$, \\ 
{\rm (iii)} $(X, X_1, X_2)$ is $(s, \Sigma, \sigma)$-stable. 
\end{proposition}

We refer to \cite{vM} for related topics in infinite-dimensional topology. 

\subsection{The spaces of embeddings into 2-manifolds} \mbox{} \vskip 2mm
First we summarize the stability property and the class property of embedding spaces. Suppose $(X, d)$ and $(Y, \rho)$ are metric spaces. An embedding $f : X \to Y$ is said to be $L$-Lipschitz ($L \geq 1$) if $\frac{1}{L} d(x, y) \leq \rho(f(x), f(y)) \leq L d(x, y)$ for any $x, y \in X$. 

\begin{lemma} $($\cite[Theorems 1.2]{SW}$)$ Suppose $M$ is a Euclidean PL $2$-manifold and $K \subset X$ are compact subpolyhedra of $M$. If $\dim (X \setminus K) \geq 1$, then the triples $({\mathcal E}_K(X, M), {\mathcal E}_K^{{\rm LIP}}(X, M), {\mathcal E}_K^{{\rm PL}}(X, M))$ and $({\mathcal E}_K(X, M)^{\ast}, {\mathcal E}_K^{{\rm LIP}}(X, M)^{\ast}, {\mathcal E}_K^{{\rm PL}}(X, M)^{\ast})$ are $(s, \Sigma, \sigma)$-stable. 
\end{lemma}

\begin{lemma}
{\rm (1)} Suppose $X$ is a compact metric space, $K$ is a closed subset of $X$ and $Y$ is a locally compact, separable metric space. Then {\rm (i)} ${\mathcal E}_K(X, Y)$ is separable, completely metrizable, and {\rm (ii)} ${\mathcal E}^{{\rm LIP}}_K(X, Y)$ is $\sigma$-compact. \\ 
{\rm (2)} $($\cite{Ge}$)$ If $X$ is a compact polyhedron, $K$ is a subpolyhedron of $X$, and $Y$ is a locally compact polyhedron, then ${\mathcal E}^{{\rm PL}}_K(X, Y)$ is $\sigma$-fd-compact.
\end{lemma}

\begin{proof}
(1) (i) ${\mathcal C}(X, Y)$ is completely metrizable by the sup-metric, and ${\mathcal E}(X, Y)$ is $G_{\delta}$ in ${\mathcal C}(X, Y)$. 

(ii) For $L \geq 1$ let ${\mathcal E}^{{\rm LIP}(L)}(X, Y)$ denote the subspace of $L$-Lipschitz embeddings. If we write $Y = \cup_{n = 1}^{\infty} Y_n$ ($Y_n$ is compact and $Y_n \subset {\rm Int} \, Y_{n+1}, n \geq 1$), then ${\mathcal E}^{{\rm LIP}}(X, Y) = \cup_{n =1}^{\infty} {\mathcal E}^{{\rm LIP}(n)}(X, Y_n)$. Since ${\mathcal E}^{{\rm LIP}(n)}(X, Y_n)$ is equicontinuous and closed in ${\mathcal C}(X, Y_n)$, it is compact by Arzela-Ascoli Theorem (\cite[Ch. XII. Theorem 6.4]{Du}). Hence ${\mathcal E}^{{\rm LIP}}(X, Y)$ is $\sigma$-compact. 
\end{proof}

For the proper PL-embedding case we need some basic facts:

\begin{fact}
(1) Suppose $A$ is a PL disk (or a PL arc) and $a \in {\rm Int} \, A$. Then there exists a map $\phi : {\rm Int} \, A \to
{\cal H}^{\rm PL}_{\partial A}(A)$ such that $\phi_x(a) = x
\ (x \in {\rm Int} \, A)$ and $\phi_a = id_A$. \\ 
(2) Suppose $N$ is a PL 1-manifold with $\partial N = \emptyset$, $Y$ is a compact subpolyhedron of $N$, $U$ is an open neighborhood of $Y$ in $N$. Then there exists an open neighborhood ${\cal U}$ of $i_Y$ in ${\cal E}^{\rm PL}(Y,
N)$ and a map $\phi : {\cal U} \to {\cal H}^{\rm PL}_{N \setminus U}(N)$ such that $\phi(f)|_Y = f$ and $\phi(i_Y) = id_N$. \\ 
(3) Suppose $M$ is a PL 2-manifold, $N$ is a compact 1-submanifold of $\partial M$ and $U$ is an open neighborhood of $N$ in $M$. Then
there exists an open neighborhood ${\cal U}$ of $id_{\partial M}$ in ${\cal H}^{\rm PL}_{\partial M \setminus N}(\partial
M)$ and a map $\phi : {\cal U} \to {\cal H}_{M \setminus U}^{\rm PL}(M)$ such that $\phi(f)|_{\partial M} = f$ and $\phi(id_{\partial M}) = id_M$. \\ 
(4) Suppose $M$ is a PL 2-manifold, $Y$ is a compact subpolyhedron of $\partial M$
and $U$ is an open neighborhood of $Y$ in $M$. Then there exists an open neighborhood ${\cal V}$ of $i_Y$ in ${\cal
E}^{\rm PL}(Y, \partial M)$ and a map $\phi : {\cal V}
\to {\cal H}^{\rm PL}_{M \setminus U}(M)$ such that $\phi(g)|_Y = g$ and $\phi(i_Y) = id_M$. 
\end{fact}

\noindent {\it Comment.} \ (3) Using a PL-collar of $\partial M$ in $M$, the assertion follows from the following remarks: \\
(3-i) If $A$ is a PL arc (or a PL open arc), then there exists a map $\phi : {\cal H}^{\rm PL}_+(A) \to {\cal H}^{\rm PL}(A \times [0,1])$ such that $\phi(f)$ is an isotopy from $f$ to $id_A$ (i.e. $\phi(f)(x, t) =  (\ast, t)$, $\phi(f)(x, 0) = f(x)$ and $\phi(f)(x, 1) = (x, 1)$) for each $f \in {\cal H}^{\rm PL}_+(A)$ and $\phi(id_A) = id_{A \times [0, 1]}$. \\ 
(3-ii) Suppose $S$ is a PL circle. Then there exists an open neighborhood ${\cal U}$ of $id_S$ in ${\cal H}^{\rm PL}(S)$ and a map $\phi : {\cal U} \to {\cal H}^{\rm PL}(S \times [0, 1])$ such that $\phi(f)$ is an isotopy from $f$ to $id_S$ for each $f \in {\cal U}$ and $\phi(id_S) = id_{S \times [0, 1]}$. 

In (3-i) we may assume that $A = [0, 1]$ (or $A = {\Bbb R}$). Then $\phi(f)$ is defined as the linear isotopy $\phi(f)(x, t) =
((1-t)f(x) + tx, t)$.  

(4) This follows from (2) and (3). \qed 

\begin{lemma}
If $M$ is a PL $2$-manifold and $K \subset X$ are compact subpolyhedra of $M$, then 
{\rm (i)} ${\mathcal E}_K(X, M)^\ast$ is  completely metrizable and 
{\rm (ii)} ${\mathcal E}^{{\rm PL}}_K(X, M)^\ast$ is $\sigma$-fd-compact.
\end{lemma}

\begin{proof}
(i) ${\mathcal E}_K(X, M)^\ast$ is $G_{\delta}$ in ${\mathcal E}_K(X, M)$. \\
(ii) We may assume that $K = \emptyset$. It suffices to show that each $f \in {\cal E}^{\rm PL}(X, M)^\ast$ has a $\sigma$-fd-compact neighborhood. Since ${\cal E}_K^{\rm PL}(X, M)^\ast \cong {\cal E}_K^{\rm PL}(f(X), M)^\ast$, we may assume that $f = i_X$. Choose a sequence of small collars $C_n$ of $\partial M$ in $M$ pinched at $Y = X \cap \partial M$ such that $C_n$ becomes thinner and thinner and also the angle between ${\rm Fr}_M C_n$ and $\partial M$ at ${\rm Fr}_{\partial M} Y$ becomes smaller and smaller as $n \to \infty$.
Let $M_n = cl(M \setminus C_n)$. Then ${\cal E}^{\rm PL}_Y(X, M)^\ast = \cup_n \, {\cal E}^{\rm PL}_Y(X, M_n)$ and ${\cal E}^{\rm PL}(Y, \partial M)$ are $\sigma$-fd-compact by \cite{Ge}. 

By Fact 4.2.(4) there exists an open neighborhood ${\cal V}$ of $i_Y$ in ${\cal E}^{\rm PL}(Y, \partial M)$ and a map $\phi : {\cal V} \to {\cal
H}^{\rm PL}(M)$ such that $\phi(g)|_Y = g$ and $\phi(i_Y) = id_M$. Let $\psi : {\cal E}^{\rm PL}(X, M)^\ast \to {\cal E}^{\rm
PL}(Y, \partial M)$ be the restriction map, $\psi(f) = f|_Y$ and let ${\cal U} = \psi^{-1}({\cal V})$. Then $\Phi : {\cal V} \times {\cal E}^{\rm PL}_Y(X, M)^\ast \to {\cal U}$, $\Phi(g, h) = \phi(g)h$, is a homeomorphism with the inverse
$\Phi^{-1}(f) = (f|_Y, \phi(f|_Y)^{-1}f)$. Hence ${\cal U}$ is also $\sigma$-fd-compact. This implies the conclusion.  
\end{proof}

Next we verify the ANR-condition and the h.n.\,complement condition. 

\begin{fact} (\cite{Ga}, \cite{GH}) 
Suppose $M$ is a compact PL 2-manifold and $X$ is a compact subpolyhedron of $M$. 
Then ${\cal H}_X^{{\rm PL}}(M)$ has the h.n.\,complement in ${\cal H}_X(M)$. 
\end{fact}

\noindent {\it Comment.} By \cite[p10]{Ga} (a comment on a relative version) ${\cal H}_X^{{\rm PL}}(M)$ is (uniformly) locally contractible. Since  ${\cal H}_X(M)$ is an ANR, by \cite{GH} ${\cal H}_X^{{\rm PL}}(M)$ has the h.n.\,complement in ${\cal H}_X(M)$. Note that  in dimension 2, the local contractibility of ${\cal H}_X^{\rm PL}(M)$ at $id_M$ simply reduces to the case where $X = \emptyset$ by the following splitting argument: \\
(1) We may assume that $X$ has no isolated points in ${\rm Int}\, M$. If $X$ has the isolated points $x_i \ (i = 1, \cdots, n)$ in ${\rm Int}\, M$, then we can choose mutually disjoint PL disk neighborhood $D_i$ of $x_i$ in ${\rm Int}\, M \setminus X_0$, where $X_0 = X \setminus \{ x_1, \cdots, x_n \}$. By Fact 4.2.(1) there exists a map $\phi :
\prod_{i=1}^n {\rm Int} D_i \to {\cal H}^{\rm PL}_{X_0}(M)$ such that $\phi(y_1, \cdots, y_n)(x_i) = y_i$ and $\phi(x_1, \cdots, x_n) =
id_M$. Then ${\cal U} = \{ f \in {\cal H}^{\rm PL}_{X_0}(M) \, : \, f(x_i) \in {\rm Int} \, D_i \ (i = 1, \cdots, n) \}$ is an open neighborhood of $id_M$ in ${\cal H}^{\rm PL}_{X_0}(M)$ and $\Phi : (\prod {\rm Int} D_i) \times {\cal H}^{\rm PL}_X(M) \to {\cal U}, \ \Phi(y_1, \cdots, y_n, g) = \phi(y_1, \cdots, y_n)g$, is a homeomorphism with the inverse $\Phi^{-1}(f) = (f(x_1), \cdots, f(x_n), \\ \phi(f(x_1), \cdots, f(x_n))^{-1}f)$. Hence if ${\cal H}_{X_0}^{{\rm PL}}(M)$ is locally contractible, then ${\cal H}_X^{{\rm PL}}(M)$ is also locally contractible. \\
(2) Cutting $M$ along ${\rm Fr}_M X$ we may assume that $X \subset \partial M$. \\
(3) By Fact 4.2.(4) there exists an open neighborhood ${\cal V}$ of $i_X$ in ${\cal E}^{\rm PL}(X, \partial M)$ and a map $\phi : {\cal V} \to {\cal H}^{\rm PL}(M)$ such that $\phi(g)|_X = g$ and $\phi(i_X) = id_M$. Let $\psi : {\cal H}^{\rm PL}(M) \to {\cal E}^{\rm PL}(X, \partial M)$ be the restriction map, $\psi(f) = f|_X$ and let ${\cal U} = \psi^{-1}({\cal V})$. Then ${\cal U}$ is an open neighborhood of $id_M$ in ${\cal H}^{\rm PL}(M)$ and $\Phi : {\cal V} \times {\cal H}^{\rm PL}_X(M) \to {\cal U}$, $\Phi(g, h) = \phi(g)h$, is a homeomorphism with the inverse $\Phi^{-1}(f) = (f|_X, \phi(f|_X)^{-1}f)$. Since ${\cal H}^{\rm PL}(M)$ is locally contractible \cite{Ga}, ${\cal H}^{\rm PL}_X(M)$ is also locally contractible.  \qed
\vskip 2mm 

Suppose $M$ is a PL $2$-manifold and $K \subset X$ are compact subpolyhedra of $M$. 

\begin{lemma}
{\rm (1) (i)} ${\mathcal E}_K(X, M)^{\ast}$ is an ANR and {\rm (ii)} ${\mathcal E}_K^{{\rm PL}}(X, M)^{\ast}$ has the h.n.\ complement in ${\mathcal E}_K(X, M)^{\ast}$. \\ 
{\rm (2)} {\rm (i)} ${\mathcal E}_K(X, M)$ is an ANR and {\rm (ii)} ${\mathcal E}_K^{{\rm PL}}(X, M)$ has the h.n.\ complement in ${\mathcal E}_K(X, M)$. 
\end{lemma}

\begin{proof}
(1)(i) For every $f \in {\mathcal E}_K(X, M)^\ast$, take a compact PL 2-submanifold neighborhood $N$ of $f(X)$ in $M$ and consider the map $\pi : {\mathcal H}_{K \cup (M \setminus {\rm Int}_M N)}(M) \to {\mathcal E}_K(X, M)^\ast$, $\pi(h) = hf$. By Theorem 1.1 there exists an open neighborhood ${\cal U}$ of $f$ in ${\mathcal E}_K(X, M)^\ast$ and a map $\varphi : {\cal U} \to {\mathcal H}_{K \cup (M \setminus {\rm Int}_M N)}(M)$ such that $\pi\phi(g) = g \ (g \in {\cal U})$. Since ${\mathcal H}_{K \cup (M \setminus {\rm Int}_M N)}(M) \cong {\mathcal H}_{K \cup {\rm Fr}_M N}(N)$ is an ANR by Lemma 3.2, so is ${\cal U}$. Hence ${\mathcal E}_K(X, M)^\ast$ is an ANR.

(ii) By Fact 4.1.(i) it suffices to show that every $f \in {\cal E}_K(X, M)^\ast$ admits a neighborhood ${\cal U}$ and a homotopy $F_t : {\cal U} \to {\cal E}_K(X, M)^\ast$ such that $F_0 = i_{\cal U}$ and $F_t({\cal U}) \subset {\mathcal E}_K^{{\rm PL}}(X, M)^{\ast} \ (0 < t \leq 1)$. Take a compact PL 2-submanifold $N$ of $M$ with $f(X) \subset U \equiv {\rm Int}_M N$. Let $\phi : {\cal U} \to {\cal H}_{K \cup (M \setminus U)}(M)$ be given by Theorem 1.1. Since $({\cal H}_{K \cup (M \setminus U)}(M), {\cal H}_{K \cup (M \setminus U)}^{{\rm PL}}(M)) \cong ({\cal H}_{K \cup (N \setminus U)}(N), {\cal H}_{K \cup (N \setminus U)}^{{\rm PL}}(N))$, by Fact 4.2 we have an absorbing homotopy $\chi_t$ of ${\cal H}_{K \cup (M \setminus U)}(M)$ into ${\cal H}_{K \cup (M \setminus U)}^{{\rm PL}}(M)$. There exists a $h \in {\cal H}_{K \cup (M \setminus U)}(M)$ such that $hf \in {\cal E}^{{\rm PL}}_K(X, M)^\ast$. Define $F_t$ by $F_t(g) = \chi_t(\phi(g)h^{-1})h f \ (g \in {\cal U})$. 

(2) There exists an $f \in {\cal E}^{{\rm PL}}_K(X, M)$ with $f(X \setminus K) \subset {\rm Int} \, M$. It induces a homeomorphism $({\cal E}_K(f(X), M)$, ${\cal E}^{{\rm PL}}_K(f(X), M)) \cong ({\cal E}_K(X, M), {\cal E}^{{\rm PL}}_K(X, M)) \, : \, g \mapsto gf$. Hence we may assume that $X \setminus K \subset {\rm Int} \, M$. Pushing towards ${\rm Int} \, M$ using a collar of $\partial M$ pinched on $\partial M \cap K$, it follows that ${\mathcal E}_K(X, M)^\ast$ has the h.n.\ complement in ${\mathcal E}_K(X, M)$. Thus (i) follows from (1)(i) and Fact 4.1.(ii), and (ii) follows from (1)(ii), Fact 4.1.(iii) and ${\cal E}^{{\rm PL}}_K(X, M)^\ast \subset {\mathcal E}^{{\rm PL}}_K(X, M)$.
\end{proof}

Theorem 1.2 follows from Proposition 4.1 and the above lemmas. For the proper embeddings we have a pair version.

\begin{proposition}
If ${\rm dim} \, (X \setminus \, K) \geq 1$, then $({\mathcal E}_K(X, M)^{\ast}, {\mathcal E}_K^{{\rm PL}}(X, M)^{\ast})$ is an $(s, \sigma)$-manifold.
\end{proposition}

\begin{remark}  In general, ${\mathcal E}_K^{{\rm LIP}}(X, M)^{\ast}$ is {\it not} $\sigma$-compact. For example, suppose $X$ is a proper arc in $M$ and $K = \partial X$. If ${\cal E}_K^{\rm LIP}(X, M)^\ast = \cup_{i \geq 1}{\cal F}_i$, ${\cal F}_i$ is compact, then $F_i = \{ f(x) \, | \, f \in {\cal F}_i, x \in X \}$ is a compact subset of $M$ with $F_i \cap \partial M = K$. By a simple diagonal argument we can define an $f \in {\cal E}_K^{\rm LIP}(X, M)^\ast$ such that $f(X) \not\subset F_i$ for each $i \geq 1$. Figure 2 indicates how to define such an $f$ near an end point of $X$. 
\end{remark} 

\[ \fbox{Figure 2.} \]

\end{document}